\documentclass[11pt]{amsart}
\usepackage{amsmath,amscd,amssymb,amsfonts}
\theoremstyle{plain}
\newtheorem{thm}{Theorem}[section]

\newtheorem{definition}[thm]{Definition}
\theoremstyle{remark}

\newtheorem{remark}[thm]{Remark}

\numberwithin{thm}{section}
\numberwithin{equation}{section}

\newcommand{\sD}{{\mathcal D}}

\newcommand{\C}{{\mathbb C}}

\renewcommand{\P}{{\mathbb P}}
\newcommand{\Q}{{\mathbb Q}}

\newcommand{\Z}{{\mathbb Z}}

\begin{document}

\title[Differential Equations associated to Families of Algebraic Cycles]
{Differential Equations associated to Families of Algebraic Cycles}
 
\author{Pedro Luis del Angel}
\address{Pedro Luis del Angel, CIMAT, Guanajuato, Mexico}
\email{luis@cimat.mx}

\author{Stefan M\"uller-Stach}
\address{Stefan M\"uller-Stach, Institut f\"ur Mathematik, Fachbereich 08, 
Johannes Gutenberg--Universit\"at Mainz, Deutschland}
\email{mueller-stach@uni-mainz.de}

\date{\today}
\dedicatory{Dedicated to Sevin Recillas}

\begin{abstract} 
We develop a theory of differential equations associated to families
of algebraic cycles in higher Chow groups (i.e., motivic cohomology groups). 
This formalism is related to inhomogenous Picard--Fuchs type differential equations.
For families of K3 surfaces the corresponding non--linear ODE turns out to be similar to 
Chazy's equation.
\end{abstract}
\subjclass{ 14C25, 19E20 }
\thanks{We are grateful to CONACYT (grant 37557-E), CIMAT
DFG Heisenberg and Schwerpunkt program, the Fields Institute,  
McMaster University, Universit\"at Mainz 
and Universit\"at Essen for supporting this project}
\maketitle
\section{Introduction}  

Around 1900 R. Fuchs~\cite{F} discovered a connection between 
non--linear second order ODE of type Painlev\'e VI~\cite{P} and integrals of
holomorphic forms over non--closed paths on  the Legendre family of elliptic curves. 
During the whole 20th century the Painlev\'e VI equation has played a prominent role
in mathematics and physics, see \cite{umemura}. About 100 years later,     
Y.I.Manin~\cite{M} found a framework in which inhomogenous   
Picard--Fuchs $\mu$--equations and non--linear equations of type 
Painlev\'e VI can be connected to mathematical physics and 
the theory of integrable systems. Inspired by his work and the 
earlier work of Griffiths~\cite{G4} and Stiller~\cite{St}
about differential equations satisfied by normal functions, 
the authors~\cite{dAMS} have looked at inhomogenous equations 
in the case of the higher Chow group $CH^2(X,1)$ of K3 surfaces. \\

In this paper we study differential equations arising from families of 
algebraic cycles in higher Chow groups $CH^p(Y,n)$ of projective manifolds $Y$ 
with $2p-n-1=\dim(Y)=d$. Our goal is to develop a theory of 
differential equations associated to each family ${\mathcal Z}/B$ of cycles in higher Chow groups over a 
quasi--projective base variety $B$. In \cite{dAMS} we suggested 
to use any Picard--Fuchs operator ${\mathcal D}$ of the local system underlying the smooth family 
$$
f: X \longrightarrow B
$$
and the new invariants are given by the assignment
$$
{\mathcal Z}/B \mapsto g(t):={\mathcal D} \nu_{{\mathcal Z}/B}(\omega),
$$
where $\nu_{{\mathcal Z}/B}$ is the \emph{normal function} associated to the 
family ${\mathcal Z}/B$ and $\omega$ a relative smooth $d$--form. In particular the function $g(t)$ depends on $\omega$ and 
the choice of Picard--Fuchs operator ${\mathcal D}$. This construction can be used in the following way: 
If we want to prove that a family of cycles is non--trivial, i.e., its Abel--Jacobi image is non--zero modulo torsion 
then it is sufficient to show that $g(t)={\mathcal D} \nu_{{\mathcal Z}/B}(\omega)$ is not zero for some choice 
of $\omega$ \cite{dAMS}.
 
In section~\ref{diffeq} we discuss the differential equations satisfied by admissible normal functions using Picard--Fuchs operators. 
This gives us the possibility to investigate the relation between the field of definition of ${\mathcal Z}$ and the coefficients of $g$ in 
section~\ref{cycles}. 
We restrict ourselves to the case of varieties with trivial canonical bundle, 
where the choice of $\omega$ is unique up to an invertible function on the base. However this restriction 
is not necessary, in the general case we will obtain a vector valued invariant. We prove under these assumptions:

\begin{thm} If, under these assumptions, ${\mathcal Z}$ and $\omega_{X/B}$ are defined over an algebraically 
closed field $K \supset \overline{\Q}$, then ${\mathcal D}$ and $g(t)$ have coefficients in $K$ and $g(t)$ is an algebraic function of $t$. 
\end{thm}

The differential equation 
$$
{\mathcal D} \nu_{{\mathcal Z}/B}(\omega)=g(t)
$$
for $\nu_{{\mathcal Z}/B}$ thus contains in general some interesting information 
about the cycle ${\mathcal Z}$, provided that the monodromy and the cycle under consideration is non--trivial.
In particular if the set of singularities (i.e., poles and algebraic branch points) of $g$ are fixed then there is only a 
countable set of possibilities for the coefficients.

In section~\ref{examples} we recall the case of dimension $1$, where this inhomogenous equation is related to 
the Painlev\'e VI equation, a second order ODE having the Painlev\'e property, i.e., no movable branch points
and essential singularities. In dimension $2$ the inhomogenous equation is of the form
$$
{\mathcal D} \int_{a(\lambda(t))}^{b(\lambda(t))} dx 
\int_{c(x,\lambda(t))}^{d(x,\lambda(t))} F(x,y,\lambda(t)) dy =g(t), 
$$
with algebraic functions $a,b,c,d$. The resulting non--linear ODE is -- after some substitutions -- of the form 
$$
\lambda '''(t) = A(\lambda)+B(\lambda) \lambda '\lambda '' + C(\lambda)(\lambda ')^3. 
$$
This equation is a variant of Chazy's equation~\cite[page 319]{C}, a third order ODE with the Painlev\'e property.  
In the study of isomonodromic deformations such PDE also arise, see \cite{boalch}. 
In future work we will come back to the Painlev\'e property in our setup.

\section{Cycle class maps from higher Chow groups to Deligne cohomology}
\label{cycleclass}

Higher Chow groups \cite{B} can be defined using the algebraic $n$--cube
$$
\square^n=({\mathbb P}^1_F \setminus \{1\})^n. 
$$
The $n$--cube has $2^n$ codimension one faces, defined by $x_i=0$ and 
$x_i=\infty$, for $1 \le i \le n$,
and the boundary maps are given by 
$$
\partial=\sum_{i=1}^n (-1)^{i-1} (\partial_i^0 -\partial_i^\infty),
$$ where $\partial_i^0$ and $\partial_i^\infty$
denote the restriction maps to the faces $x_i=0$ and $x_i=\infty$. 
Then $Z_c^p(X,n)$ is defined to be the quotient of the group
of admissible cycles in $X \times \square^n$ by the group of degenerate cycles, see \cite{B}.
We use the notation $CH^p(X,n)$ for the $n$--homology of the complex $Z_c^p(X,\cdot)$. 
There are cycle class maps 
$$
c_{p,n}: CH^p(X,n) \longrightarrow H_{\mathcal D}^{2p-n}(X,{\mathbb Z}(p))
$$
constructed by Bloch in \cite{B3} using Deligne cohomology \cite{CMP} with supports
and a spectral sequence construction. They can be realized explicitly by 
Abel--Jacobi type integrals if $X$ is a complex, projective manifold~\cite{KLM}.
If we restrict to cycles homologous to zero, then we obtain Abel--Jacobi type maps
$$
c_{p,n}: CH^p(X,n)_{\rm hom} \longrightarrow J^{p,n}(X)=\frac{H^{2p-n-1}(X,\C)}{F^p + H^{2p-n-1}(X,\Z)},
$$
where $J^{p,n}(X)$ are \emph{generalized intermediate Jacobians} \cite{KLM}.
These are complex manifolds and vary holomorphically in families like Griffiths' intermediate Jacobians \cite{CMP}.

\section{Differential Equations associated to families of algebraic cycles}
\label{diffeq}

In this section we study \emph{differential equations} arising from families of 
algebraic cycles. Assume that we are in the following setup:

Let $f:X \to B$ a smooth, projective family of manifolds with trivial canonical bundle (e.g. Calabi--Yau) of 
relative dimension $d=2p-n-1$ over a smooth, quasi--projective curve $B$ 
with compactification $\overline B$. Assume that $f$ is defined over an algebraically closed field $K \subseteq \C$.

We fix a base point $o \in B$ and a local parameter
$t$ around $o$ with  $t \in K(B)$, the function field of $B$, so that $dt$ is a basis of 
$\Omega_{B,o}^1$ and $\frac{\partial}{\partial t}$ the corresponding vector field. 
Let ${\mathbb H}$ be the local system associated to the primitive part of $R^{d}f_*{\mathbb C}$. Its stalks 
consist of cohomology groups $H^{d}_{\rm pr}(X_t,{\mathbb C})$ for $t \in B$. 
We assume that ${\mathbb H}$ has an irreducible monodromy representation with unipotent local behaviour around each 
point at infinity. Denote by ${\mathcal H}$ the holomorphic vector bundle with sheaf of sections
${\mathcal H}= {\mathbb H} \otimes {\mathcal O}_B$ and Gau{\ss}--Manin connection $\nabla$. The Hodge pairing is denoted by  
$\langle -,- \rangle: {\mathcal H} \otimes  {\mathcal H} \to {\mathcal O}_B$. Together with the Hodge filtration $F^\bullet$ 
this data defines a polarized VHS on $B$. 

We choose a non--zero holomorphic section $\omega \in H^0(B,F^d {\mathcal H})$ and denote by 
$$
{\mathcal D}_{\rm PF}=\frac{d^m}{dt^m}+a_{m-1}(t)\frac{d^{m-1}}{dt^{m-1}}+\cdots +a_0(t)
$$ 
the Picard--Fuchs operator corresponding to ${\mathbb H}$ in the local basis 
$\omega$, $\nabla_t \omega$, ..., $\nabla^m_t \omega$ with rational functions $a_i(t)$.

Assume furthermore that we have a cycle ${\mathcal Z} \in CH^p(X,n)$ such
that each restriction $Z_t:={\mathcal Z}|_{X_t} \in CH^p_{\rm hom}(X_t,n)$ is a well--defined
cycle, in other words we have a single--valued family of algebraic cycles over $B$. 
This implies that we have a well--defined normal function
$$
\nu \in H^0(B,{\mathcal J}^{p,n}), \quad \nu(t):=c_{p,n}(Z_t),
$$
i.e., a holomorphic cross section of the bundle ${\mathcal J}^{p,n}$ 
of generalized intermediate Jacobians. Locally on $B$ near the point $o$ in the analytic topology we may choose 
a lifting $\tilde \nu$ of $\nu$ as a holomorphic cross section of ${\mathcal H}/F^p$ or of ${\mathcal H}$ using 
identical notation. 

Any cycle $Z_t \in CH_{\rm hom}^p(X_t,n)$ defines a extension of two pure Hodge structures~\cite{KLM}
$$
0\to H^{d}_{\rm pr}(X_t) \to E_t \to {\mathbb Z}(-p) \to 0.
$$
For each $t$ the extension class of this sequence in the category of mixed Hodge structures is the Abel--Jacobi map of $Z_t$ in 
$J^{p,n}(X_t)$ \cite{KLM}. For varying $t$, $E_t$ defines a local system ${\mathbb E}$ over $B$ 
which is an extension of ${\mathbb H}$  by a trivial local system of rank one. ${\mathcal E} ={\mathbb E} \otimes {\mathcal O}_B$ carries a holomorphic flat connection $\tilde \nabla$ extending $\nabla$ and a filtration $F^\bullet$ by subbundles extending the one on ${\mathcal H}$.  
Let $\hat {\mathcal E}$ be the Deligne extension \cite{D} of ${\mathcal E}$ to $\bar B$. 
For technical reasons we will assume that $\nu $ is {\sl admissible}, i.e., the extension of MHS above is admissible in the sense of M. Saito~\cite{S1}. 
This means in particular (see loc. cit.):
\begin{itemize}
\item The Hodge filtration $F^\bullet$ on ${\mathcal E}$ extends to the Deligne extension $\hat {\mathcal E}$ with locally free graded quotients,
\item The relative monodromy weight filtration extends. 
\end{itemize}
We will use the first property in an essential way, which implies that $\nu$ has moderate growth at infinity as we will see in the proof. 
The admissibility condition is always satisfied in the geometric case when $n=0$ by Steenbrink and Zucker 
\cite[sect. 3]{SZ}. In general for $n\ge 1$ it is not well--understood. However extendable normal functions in the sense of 
M. Saito are admissible by \cite[Prop. 2.4]{S1}.

Since $E_t$ is an extension by a pure Hodge structure ${\mathbb Z}(-p)$ of type $(p,p)$, we have
$$
{\mathcal E}/F^p = {\mathcal H}/F^p.
$$ 
After further lifting, we can view $\tilde \nu$ by abuse of notation as a multivalued holomorphic section of either ${\mathcal E}$ or 
${\mathcal H}$. In each case it is well--defined modulo $F^p$ only.
It is not a flat section for $\nabla$ unless the cycle has trivial Abel--Jacobi invariant. 

\begin{definition} The truncated normal function $\bar \nu$ is defined as $\bar \nu:= \langle \tilde \nu,\omega \rangle$. 
\end{definition}

Formulas for $\bar \nu$ are given in \cite{KLM} using so--called {\sl membrane integrals}.
Note that $\bar \nu$ does not depend on the lifting $\tilde \nu$ if $p \ge 1$, since the holomorphic $d$--form 
$\omega$ has only non--zero Hodge pairing with $(0,d)$--classes which are never contained in $F^p$. 

\begin{thm} \label{moderategrowth} Let $\nu$ be an admissible higher normal function as above. Then
$\bar \nu$ is a multivalued function on $B$. Furthermore we have ${\mathcal D}_{\rm PF}\bar \nu=g(t)$ 
for some single--valued holomorphic function $g(t)$ on $B$. The Picard--Fuchs equation for $\bar \nu$ is homogenous:
$$
\left(\frac{d}{dt} - \frac{g'(t)}{g(t)}\right) \cdot {\mathcal D}_{\rm PF}\bar \nu =0.
$$ 
In particular, the holomorphic function $g(t)$ extends to a rational function on $\overline B$.
\end{thm} 

\proof First we show that $\tilde \nu$ can be chosen flat when considered as a section of $({\mathcal E},\tilde \nabla)$.
We use Carlson's extension theory of Hodge structures which in our case says that the extension class of the sequence 
$$
0\to H^{d}_{\rm pr}(X_t) \to E_t \to {\mathbb Z}(-p) \to 0
$$
in the category MHS is given (up to a sign) by an integral lifting $s_{\Z}$ of $1 \in {\mathbb Z}(-p)$~\cite{KLM}. 
The Abel--Jacobi invariant is then obtained by projecting $s_{\Z}$ into 
$$
J^{p,n}(X_t)=\frac{H^{d}_{\rm pr}(X_t,\C)}{F^p+H^{d}_{\rm pr}(X_t,\Z)}=\frac{E_t(\C)}{F^p+E_{t}(\Z)}.
$$ 
We use that $\tilde \nu$ is defined as a current of integration defined in 
\cite{KLM}. In the classical situation, i.e., $n=0$ it is given by the current $\alpha \mapsto \int_{\Gamma_t} \alpha$, 
which is dual to a relative homology class of $\Gamma_t$ in $H_d(X_t,|Z_t|,{\mathbb Z})$. More formally, one has the long exact sequence 
$$
H^{2p-1}(X_t \setminus |Z_t|,{\mathbb Z}) \to H^{2p}_{|Z_t|}(X_t,{\mathbb Z}) \to H^{2p}(X_t,{\mathbb Z}).
$$
Since $Z_t$ is homologous to zero, its fundamental class in $H^{2p}_{|Z_t|}(X_t,{\mathbb Z})$ can be non--uniquely lifted to 
a class $s_{\Z}$ in $H^{2p-1}(X_t \setminus |Z_t|,{\mathbb Z})$. $s_{\Z}$ is unique up to elements in $H^{2p-1}(X_t,{\mathbb Z})$
which however vanish in the intermediate Jacobian and represents therefore $\tilde \nu$ up to the choices in $F^p+H^{2p-1}(X_t,\Z)$ 
by Carlson's theory. This proves the assertion, since integral classes are always flat. 

In the case $n \ge 1$ the argument is similar. The support $|Z_t|$ is a subset of $X \times \Box^n$.
The mixed Hodge structure associated to ${\mathbb E}_t$ is a subquotient of the relative cohomology group 
$H^{2p-1}(U_t,\partial U_t)$, where $U_t:=X \times \Box^n \setminus |Z_t|$ and $\partial U_t:= U_t \cap \partial \Box^n$. 
One has then an exact sequence with integral coefficients \cite[(6.1)]{KLM} 
$$
0 \to H^{d}(X_t) \to H^{2p-1}(U_t,\partial U_t) \to {\rm ker}(\beta) \to H^{2p-n}(X_t),
$$
where $\beta$ is the map 
$$
\beta: H^{2p}_{|Z_t|}(X_t \times \Box^n)^\circ \to H^{2p}_{\partial |Z_t|}(X_t \times \partial \Box^n)^\circ.
$$
The symbol $\circ$ stands for the kernel of the map forgetting supports. For any $[Z_t] \in {\rm ker}(\beta)$ we obtain an extension
$$
0\to H^{d}_{\rm pr}(X_t) \to E_t \to {\mathbb Z}(-p) \to 0
$$
as a subquotient. As in the case $n=0$ we conclude that we can lift the fundamental class $[Z_t]$ to an integral class $s_\Z$ 
which coincides with $\tilde \nu$ up to the choices in $F^p+H^{d}(X_t,\Z)$. Flatness follows again from integrality.

Now, since $\tilde \nu$ becomes flat as a section of $({\mathcal E},\tilde \nabla)$, this then 
implies that $\bar \nu$ is a multi--valued solution of the homogenous Picard--Fuchs equation associated to 
$({\mathcal E},\tilde \nabla)$~\cite[Prop 8.1.4]{K}. 
Since $\bar \nu$ satisfies the inhomogenous Picard--Fuchs equation ${\mathcal D}_{\rm PF}=g(t)$, it is a solution of
$$
\left(\frac{d}{dt} - \frac{g'(t)}{g(t)}\right) \cdot {\mathcal D}_{\rm PF}\bar \nu =0.
$$ 
Therefore this equation must be the homogenous Picard--Fuchs equation associated to $({\mathcal E},\tilde \nabla)$. 
Since ${\mathcal E}$ is of algebraic origin and admissible, we deduce in addition that $g(t)$ is a rational function of 
$t$ as in the proof of \cite[Prop. 3.12]{St}. \qed

\begin{remark}
The same trick also shows that any admissible normal function in our setup is a G--function in the sense of Siegel and Andre, 
see \cite{A}. 
\end{remark}

\section{Applications to algebraic cycles}
\label{cycles}

Let $f:X \to B$ be a smooth, projective family of projective manifolds with trivial canonical bundle (e.g. Calabi--Yau) 
of dimension $d=2p-n-1$ over a smooth, quasi--projective curve $B$ with projective compactification $\overline{B}$.  
As in the previous section we are given a single--valued family of cycles ${\mathcal Z_t} \in CH^p(X_t,n)$ inducing
a well--defined normal function
$$
\nu \in H^0(B,{\mathcal J}^{p,n}), \quad \nu(t):=c_{p,n}(Z_t).
$$
We also use the same notations for the irreducible local system ${\mathbb H}$ 
of primitive cohomology and assume that its has unipotent local monodromies at infinity.
Let $\sD_{\rm PF}$ be a Picard--Fuchs operator for this family after a choice 
of $\omega \in H^0(B,F^d {\mathcal H})$.

Assume in addition that $f: X \to B$ and the cycle ${\mathcal Z}$ 
are defined over $\overline \Q$ or -- more generally -- over any algebraically closed extension field $K \supseteq \overline \Q$. 
Such a situation can for example be achieved by spreading out a cycle on a generic fiber $X_\eta$ over the field obtained by the 
compositum of its field of definition and the function field of $\eta$. In other 
words all transcendental elements in the equations of $X_\eta$ and ${\mathcal Z}$ 
occur in the coordinates of $B$. Then there is a canonical choice of a relative holomorphic 
$d$--form $\omega$ defined over $\overline \Q$. In our case, where $B$ is a curve, such a situation 
is only possible if the transcendance degree of $K$ over $\Q$ is $\le 1$.

The following theorem leads us to expect that normal functions of cycles defined over $K$ with a fixed set of critical 
points (i.e., poles) form at most a countable set.

\begin{thm} \label{countable}
The rational function $g$ has all its coefficients in $K$ under these assumption. 
\end{thm}

\proof Since $Z$ and $X$ are defined over $K$, the cohomology class of $Z$ in $F^p H_{\rm dR}^{2p-n}(X)$ and 
the extension data of VMHS in the proof of Theorem~\ref{moderategrowth} are defined over $K$. Hence the Gau{\ss}--Manin connection and 
the Picard--Fuchs operator have coefficients in $K$ as well.  
Theorem~\ref{moderategrowth} implies that $g$ is a rational function with coefficients in $K$. \endproof

\begin{remark} Our proof can be generalized to a higher dimensional base variety $B$. Then the occurring Picard--Fuchs systems
will define in general a non--principal ideal of partial differential operators. We may then assume that the transcendance degree
of $K$ is as large as $\dim(B)$. As above we can only expect single--valuedness and unipotency after 
a finite base--change. Therefore $\overline \nu$ will in general be an algebraic function over $\overline{B}$. 
Manin's example in \cite{M} already involves a square root.
Finally we want to remark that the normal functions are not necessarily uniquely 
determined by these differential equations since there may be a non--trivial monodromy 
invariant part of the cohomology. 
\end{remark}

\section{Examples}
\label{examples} 

In this section we give examples in dimensions $1$ and $2$ and relate them to classical 
non--linear ODE. For dimension $1$, consider a section 
$$
t \mapsto \left(X(t),Y(t)\right)
$$ 
of the Legendre family, written as 
$$
y^2=x(x-1)(x-t), \quad t \in \P^1 \setminus \{0,1,\infty\}
$$ 
in affine coordinates. The corresponding inhomogenous Picard--Fuchs differential equation  
can be written as 
$$
{\mathcal D} \int_\infty^{X(t)} \frac{dx}{y}=g(t)
$$
for a rational or algebraic function $g(t)$, where
$$
{\mathcal D}=t(1-t)\frac{d^2}{dt^2}+(1-2t)\frac{d}{dt}-\frac{1}{4}.
$$
Richard Fuchs~\cite{F} looked at a 4--parameter 
set of such equations of the form 
$$
t(1-t){\mathcal D} \int_\infty^{X(t)} \frac{dx}{y}=Y(t) [\alpha +\beta\frac{t}{X(t)^2}+\gamma\frac{(t-1)}{(X(t)-1)^2}
+(\delta -\frac{1}{2})\frac{t(t-1)}{(X(t)-t)^2}]
$$
with $\alpha, \beta,\gamma, \delta \in \C$. Furthermore every solution of this equation is 
also a solution of the non--linear equation Painlev\'e VI and vice versa:  
$$
P_{VI}: \quad \frac{d^2X}{dt^2}=\frac{1}{2}\left(
\frac{1}{X}+\frac{1}{X-1}+\frac{1}{X-t}\right)
\left(\frac{dX}{dt}\right)^2    -\left(
\frac{1}{t}+\frac{1}{t-1}+\frac{1}{X-t}\right)\frac{dX}{dt}
$$
$$
+\frac{X(X-1)(X-t)}{t^2(t-1)^2}
\left[ \alpha +
\beta\frac{t}{X^2}+\gamma\frac{t-1}{(X-1)^2}+
\delta\frac{t(t-1)}{(X-t)^2}\right]. 
$$
This last equation $P_{VI}$ has the \emph{Painlev\'e property}, i.e., the absence of movable essential
singularities and branch points in the set of solutions. \\ 
For dimension $2$ this correspondence can be generalized: Consider a family of K3--surfaces $X_t$ 
over $\overline{B}=\P^1$ where the general fiber has Picard number $19$. Such families were considered in 
\cite[sect. 6.2.1]{SMS} and \cite{dAMS}. In this case the Picard--Fuchs operator has order $3$ and we assume 
that the cycles consist of two irreducible components. In example \cite[sect. 6.2.1]{SMS} the components are a line and an
elliptic curve. The truncated normal function $\bar \nu$ can then always be written as an integral 
$$
\bar \nu(t)= \int_{a(\lambda(t))}^{b(\lambda(t))} dx 
\int_{c(x,\lambda(t))}^{d(x,\lambda(t))} F(x,y,\lambda(t)) dy, 
$$
where $a,b,c,d$ are algebraic functions of two variables. Here $F(x,y,\lambda(t)) dxdy$ is the 
local expression for a chosen family of relative holomorphic $2$--forms. Assuming $\lambda(t)$ is
locally biholomorphic we can write $F$ as a function of $\lambda(t)$ instead of $t$. Using the same substitution 
for all coefficient functions of the Picard--Fuchs operator, which is of order $3$ here, 
we get a non--linear third order ODE of the form
$$
\lambda '''(t) = A(\lambda)+B(\lambda) \lambda '\lambda '' + C(\lambda)(\lambda ')^3,
$$
which is similar to Chazy's equation~\cite[page 319]{C}. Non--linear ODE/PDE having the Painlev\'e property like 
Chazy's equation are related to the work of Hitchin and Boalch~\cite{boalch}, where non--linear PDE occur in the theory 
of isomonodromic deformations.

\medskip
Examples in dimension $3$ related to string theory were worked out by Morrison and Walcher \cite{MW}.

\medskip
{\bf Acknowledgement:} We thank S. Bloch, H. Esnault, M. Green, Ph. Griffiths, Y.I. Manin, D. van Straten 
and J. Walcher for several valuable discussions, remarks or letters. Referees have made useful remarks concerning 
the presentation and the content of a previous version.

\bibliographystyle{plain}

\begin{thebibliography}{99}

\bibitem{A} {\sc Y. Andr\'e:} \textit{G--functions}, Aspects of Mathematics {\bf E13}, Viehweg Verlag (1989).
\bibitem{B} {\sc S. Bloch:} \textit{Higher Chow groups: Basic definitions and properties}, 
available on Bloch's homepage. 
\bibitem{B3} {\sc S. Bloch:} \textit{Algebraic cycles and the Beilinson 
conjectures}, Cont. Math. \textbf{58} Part I, 65--79 (1986).
\bibitem{boalch} {\sc Ph. Boalch :}  \textit{Symplectic manifolds and isomonodromic deformations}, 
Adv. Math. {\bf 163} no. 2, 137--205 (2001). 
\bibitem{C} {\sc J. Chazy:} \textit{Sur les \'equations diff\'erentielles du troisi\`eme
ordre et d'ordre sup\'erieur dont l'int\'egrale g\'en\'erale a ses points critiques
fixes}, Acta Math. {\bf 34}, 317--385 (1910).
\bibitem{dAMS} {\sc P. L. del Angel and S. M\"uller-Stach:} \textit{The transcendental
part of the regulator map for $K_1$ on a family of K3 surfaces}, Duke
Math. Journal {\bf 112} No.3, 581--598 (2002). 
\bibitem{dAMS2} {\sc P. L. del Angel and S. M\"uller-Stach:} \textit{Picard-Fuchs equations, 
Integrable Systems and higher algebraic K--theory}, Calabi-Yau varieties and mirror symmetry (Toronto, ON, 2001),  43--55,
Fields Inst. Comm. {\bf 38}, Amer. Math. Soc., Providence, RI, (2003).
\bibitem{CMP} {\sc J. Carlson, S. M\"uller-Stach and Ch. Peters:} \textit{Period mappings and period domains}, 
Cambridge Studies in Advanced Math. {\bf 85}, Cambridge Univ. Press (2003).
\bibitem{D} {\sc P. Deligne:} \textit{Equations diff\'erentielles \`a points r\'eguliers singuliers}, Springer LNM {\bf 163} (1970). 
\bibitem{F} {\sc R. Fuchs:} \textit{\"Uber lineare homogene Differentialgleichungen zweiter
Ordnung mit drei im Endlichen gelegenen wesentlich singul\"aren Stellen}, Math. Ann.
{\bf 63}, 301--321 (1907).
\bibitem{G} {\sc M. Green:} \textit{Griffiths's infinitesimal invariant and the Abel--Jacobi map}, 
Journal Diff. Geometry {\bf 29}, 545--555 (1989).
\bibitem{G3}  {\sc Ph. Griffiths:} \textit{On the periods of certain rational integrals I,II,III} 
Ann. of Math. (2) {\bf 90}, 460--495 (1969) and 461--495 (1969), Publ. Math. IHES
\textbf{38}, 125--180 (1970).
\bibitem{G4}  {\sc Ph. Griffiths:}  \textit{A theorem concerning the differential equations 
satisfied by normal functions associated to algebraic cycles},  Amer. J. Math. 
{\bf 101}, 94--131 (1979). 
\bibitem{KLM}  {\sc M. Kerr, J. Lewis and S. M\"uller--Stach:} \textit{The Abel-Jacobi map for
higher Chow groups}, Compositio Math. {\bf 142}, 374--396 (2006).
\bibitem{K}{\sc V. Kulikov:} \textit{Mixed Hodge structures and singularities}, 
Cambridge Tracts in Mathematics {\bf 132}, 1998.
\bibitem{M} {\sc Y. I. Manin:} \textit{Sixth Painlev\'e equation, universal elliptic curve and
mirror of $\P^2$}, Amer. Math. Soc. Translations (2), Vol. {\bf 186},
131--151 (1998).
\bibitem{MW} {\sc D. Morrison and J. Walcher:} \textit{D--branes and normal functions},
Preprint hep-th/07094028 (2007).
\bibitem{SMS} {\sc S. M\"uller--Stach:} \textit{Constructing indecomposable motivic 
cohomology classes on algebraic surfaces}, J. Algebraic Geom. {\bf 6},  513--543 (1997).
\bibitem{shuji} {\sc S. M\"uller--Stach and S. Saito:} \textit{On $K_1$ and $K_2$
of algebraic surfaces (Appendix by Collino)}, K--theory {\bf 30}, 37--69 (2003). 
\bibitem{P} {\sc P. Painlev\'e:} \textit{Stockholm lectures}, Hermann Paris (1897), 
publicly accessible through the Cornell Digital Library.
\bibitem{S1}  {\sc M. Saito:} \textit{Admissible normal functions}, J. Algebraic Geom. 
{\bf 5}, 235--276 (1996). 
\bibitem{SZ}  {\sc J. Steenbrink and S. Zucker:} \textit{Variation of mixed Hodge 
structure I}, Invent. Math. {\bf 80}, 489--542 (1985). 
\bibitem{St} {\sc P. Stiller:} \textit{Special values of Dirichlet series, monodromy, and the periods of automorphic forms},
Memoirs of the AMS Nr. {\bf 299}, Providence (1984).
\bibitem{umemura} {\sc U. Umemura:} \textit{100 years of the Painlev\'e equation},  
Sugaku {\bf 51} no. 4, 395--420 (1999).
\end{thebibliography}
\renewcommand\refname{References}

\end{document}